\documentclass[english]{amsart}
\usepackage[T1]{fontenc}
\usepackage[latin9]{inputenc}
\usepackage{geometry}
\usepackage{array}
\usepackage{longtable}
\usepackage{float}
\usepackage{units}
\usepackage{amsmath}
\usepackage{amsthm}
\usepackage[numbers]{natbib}
\usepackage{graphicx}

\newtheorem{thm}{Theorem}


\usepackage[english]{babel}
\begin{document}
\title{Bigeometric Calculus and  Runge Kutta Method}

\author{Mustafa Riza}
\email{mustafa.riza@emu.edu.tr}
\address{Department of Physics, Eastern Mediterranean University, Gazimagusa, North Cyprus, via Mersin 10, Turkey}
\author{Bu\u{g}\c{c}e Emina\u{g}a}
\address{Department of Industrial Engineering, Girne American University, Kyrenia, North Cyprus, via Mersin 10, Turkey}
\email{bugceeminaga@gau.edu.tr}

%

\begin{abstract}
The properties of the Bigeometric or proportional derivative are presented and discussed explicitly. Based on this derivative, the Bigeometric Taylor theorem is worked out. As an application of this calculus, the Bigeometric Runge-Kutta method is derived and is applied to academic examples, with known closed form solutions, and  a sample problem from mathematical modelling in biology. The comparison of the results of the Bigeometric Runge-Kutta method with the ordinary Runge-Kutta method shows that the Bigeometric Runge-Kutta method is at least for a particular set of initial value problems superior with respect to accuracy and computation time to the ordinary Runge-Kutta method.
\end{abstract}
\maketitle


\maketitle

\section{Introduction}

After a long period of silence in the field on Non-Newtoninan Calculus introduced by  Grossmann and Katz  \cite{GK} in 1972, the field experienced a revival with the mathematically comprehensive description of Geometric multiplicative calculus by Bashirov et al. \cite{BKO}, which initiated a kickstart of numerous publications in this field. 

Grossmann and Katz have shown that it is possible to create
infinitely many calculi independently \cite{GK}. They constructed a comprehensive family of calculi, including the
Newtonian or Leibnizian  Calculus, the Geometric-multiplicative  Calculus, the Bigeometric Calculus,
and infinitely-many other calculi. In 1972, they completed their book
Non-Newtonian Calculus \cite{GK} summarising all the findings, i.e.  nine specific
non-Newtonian calculi, the general theory of non-Newtonian Calculus,
and heuristic guides for applications. 

Geometric multiplicative and Bigeometric Calculus have been becoming more and more popular in the past decade. Various applications of these two fundamental multiplicative calculi have been proposed. Exemplarily, without claim to completeness we  want to state here some of the application areas of these two multiplicative calculi. Geometric Multiplicative Calculus based on the works on \cite{GK} and \cite{BKO} was applied to various fields as modelling finance, economics and demographics using  differential equation \cite{BMTO}; numerical approximation methods \cite{ROK,OM,MG, ozyapici2013multiplicative}; biological image analysis\cite{FA,florack12}; application on literary texts \cite{BB}. In order to circumvent the restriction of  Geometric multiplicative Calculus to positive valued functions of real variable, the Geometric multiplicative Calculus was extended to complex multiplicative Calculus. After \cite{U} proposed an extension to Geometric multiplicative calculus to complex valued function of complex variable, a comprehensive mathematical description of the multiplicative complex analysis was presented by Bashirov and Riza in \cite{BR2,bashirov2011complex,BR}.  Uzer, furthermore used multiplicative calculus to model problems in antenna theory, we would like to refer to \cite{U1,U2}. Applications of Bigeometric Calculus  can be found in the field of nonlinear dynamics by the group of Rybaczuk.  We want to refer to the works \cite{rybaczuk2001concept, MR2789815,MR2449023,A}.

Especially in the area of biology, there exist numerous mathematical models  based on differential equations being quite hard to solve using standard solution methods for ordinary differential equations. The 4th order Runge-Kutta method is widely used for the numerical solution of these differential equations. Exemplarily, we want to refer to the Modelling of Gene expression using differential equations \cite{chen1999modeling}, modelling Tumor growth  \cite{agarwal11}, or modelling bacteria growth and cancer \cite{MR1750091,MR2030852}. These type of problems is used for the modelling of relative change of the numbers of cells, genes, bacteria, and viruses. Therefore, the Geometric, as well as the bigeoemetric Calculus, build the proper framework for the solution of these problems. 

The first flavour of a multiplicative derivative was given by \cite{V}, where Volterra and Hostinsky propose a derivative for a matrix function. Later on, the multiplicative derivative, underlying  this study, was proposed two times independently in  literature. The first appearance dates back to 1972  to the work  by \cite{GK} where the Bigeometric derivative was introduced, next \cite{C1} introduces in his study in 2006 the same derivative under the name proportional derivative. In the following, we will use the terminology Bigeometric derivative for the basic underlying derivative of this study. Unfortunately, a complete description for the Bigeometric, or proportional derivative is not available. Based on the extensive study of the Geometric multiplicative Calculus given by \cite{BKO}, we will state the properties of the Bigeometric  derivative and the  Bigeometric Taylor  expansion in the framework of this Calculus explicitly, exploiting the straightforward relation  between the Geometric Multiplicative and the Bigeometric derivative.  The rules of Bigeometric differentiation will be elaborated and presented explicitly in section \ref{sec2} and Appendix \ref{bgdrproof}. It is noteworthy that, compared to the Geometric multiplicative derivative, the 
 Bigeometric derivative is scale free. From the proof of the Bigeometric Taylor theorem it is self-evident why the Bigeometric taylor expansion could not be found that easily. As an application of Bigeometric calculus, the Bigeometric-Runge-Kutta method will be derived in analogy to the Geometric Runge-Kutta method proposed by Riza and Akt\"ore \cite{RA13} for the solution of Bigeometric initial value problems will be derived on the basis of the Bigeometric Taylor theorem explicitly in section \ref{sec3}. Before the conclusion of this paper, the presented method is applied to  two types of problems in section \ref{sec4}. On the one hand the Bigeometric Runge Kutta method  and the ordinary Runge Kutta method are applied to two basic problems with known closed form solutions, in order to be able to compare the two methods with respect to error estimations and performance, on the other both methods are applied to an example for mathematical modelling on tumor growth by  \cite{agarwal11} to show the general applicability of the proposed method. In both types of examples, we could identify that the Bigeometric Runge Kutta method shows a considerable better performance compared to the ordinary Runge Kutta method for the same accuracy.

\section{Properties of the Bigeometric Derivative \&  Taylor Theorem}
\label{sec2}

In this section we will first discuss the properties of the Bigeometric derivative based on the relationship between the Bigeometric derivative with the Geometric Multiplicative derivative to derive the Bigeometric Taylor theorem from the  Geometric-multiplicative Taylor theorem stated in \cite{BKO}.

\subsection{The Bigeometric Derivative and its Properties}

The Bigeometric derivative is given as:

\begin{equation}
f^{\pi}(x)=\frac{d^{\pi}f(x)}{dx}=\lim_{h\rightarrow0} \left(\frac{f\bigl(\left(1+h\right)x\bigr)}{f(x)}\right)^{\frac{1}{h}}
\label{eq:vderdef}
\end{equation}

Calculating the limit gives the  relation between the Bigeometric derivative and the ordinary derivative.

\begin{equation}
f^{\pi}(x)=\exp\left\{ x \,\frac{f'(x)}{f(x)} \right\} = \exp \left\{x  \left(\ln \circ f\right)' (x)\right\}
\label{eq:vderdef2}
\end{equation}

The complete differentiation rules of Geometric-multiplicative differentiation are presented in \cite{BKO}. As mentioned above, not all differentiation rules  of the Bigeometric derivative are presented in \cite{GR,V,C1}, therefore for sake of completeness we will state all properties of the Bigeometric derivative in the following.

 {\bf Bigeometric differentiation rules:}
 
\begin{description}
\item[Constant multiple rule]
$$\displaystyle (cf)^{\pi}(x)=(f)^{\pi}(x)$$
\item[Product Rule]
$$\displaystyle (fg)^{\pi}(x)=f^{\pi}(x)g^{\pi}(x)$$
\end{description}

\begin{description}
\item[Quotient Rule]
$$\displaystyle \left(\nicefrac{f}{g}\right)^{\pi}(x)=\nicefrac{f^{\pi}(x)}{g^{\pi}(x)}$$
\item [Power Rule]
$$\displaystyle (f^{h})^{\pi}(x)=f^{\pi}(x)^{h(x)}f(x)^{x \cdot h^{'}(x)}$$
\item[Chain Rule] 
$$\displaystyle  \left(f\circ h\right)^{\pi}(x)=f^{\pi}(h(x))^{h^{'}(x)}$$
\item[Sum Rule] 
$$\displaystyle  \left(f+g\right)^{\pi}(x)=\left(f^{\pi}(x)\right)^{\nicefrac{f(x)}{f(x)+g(x)}}\left(g^{\pi}(x)\right)^{\nicefrac{g(x)}{f(x)+g(x)}}$$
\end{description}

The proofs, carried out analogously to \cite{BKO}, are given explicitly in Appendix \ref{bgdrproof}.

For the derivation of the Bigeometric Taylor theorem the chain rule of a function of two variables will be necessary and is given as
\[
f^{\pi}(y(x),z(x))=(f_{y}^{\pi}(y(x),z(x)))^{y^{'}(x)}\cdot (f_{z}^{\pi}(y(x),z(x)))^{z^{'}(x)},
\]

with $f_y^\pi(y(x),z(x))$ denoting the partial Bigeometric derivative of $f(y(x),z(x))$ with respect to $y$, and $f_z^\pi(y(x),z(x))$ denoting the partial Bigeometric derivative of $f(y(x),z(x))$ with respect to $z$ respectively.

Proof:

\begin{multline*}
e^{x\left(\ln\left[f(y(x),z(x))\right]\right)^{'}}=\exp\left\{x \, \frac{f_{y}^{'}(y(x),z(x))\cdot y^{'}(x)+f_{z}^{'}(y(x),z(x))\cdot z^{'}(x)}{f(y(x),z(x))}\right\}=\\
=\exp\left\{x \, \frac{f_{y}^{'}(y(x),z(x))}{f(y(x),z(x))}\cdot y^{'}(x)\right\} \cdot \exp\left\{x \, \frac{f_{z}^{'}(y(x),z(x))}{f(y(x),z(x))}\cdot z^{'}(x)\right\}=\\
=f_{y}^{\pi}(y(x),z(x))^{y^{'}(x)}\cdot f_{z}^{\pi}(y(x),z(x))^{z^{'}(x)}
\end{multline*}

\subsection{Derivation of the Bigeometric Taylor Theorem}

Unfortunately, the Bigeometric Taylor theorem is not available, and the attempts by \cite{A} and \cite{ROK} show that finding the Bigeometric Taylor expansion is not straightforward.  In the following it will be clear, why finding the Bigeometric Taylor theorem is so difficult. The idea for the derivation of the Bigeometric Taylor Theorem is straightforward. We use the relation between the Bigeometric derivative and the Geometric Multiplicative derivative \eqref{eq:vderdef2}  to establish the relations for higher order derivatives as well. Once we understand the systematics of the representation of higher order Geometric-multiplicative derivatives in terms of the Bigeometric derivatives, we can substitute the higher order Geometric-multiplicative derivatives in the Geometric multiplicative Taylor theorem by a product of Bigeometric derivatives of the same function. 

First we want to find the higher order Bigeometric derivatives in terms of the Geometric Multiplicative derivative. Therefore, we will sequentially apply the relation
\[
f^\pi(x) = f^\ast(x)^x,
\]
to get the higher order derivatives.  Let us calculate the Bigeometric derivatives in terms of the Geometric Multiplicative derivatives up to order three. Using the power rule of the Geometric Multiplicative derivative \cite{BKO} in the form $\left((f^{*})^{h(x)}\right)^{*}=\left(f^{**}\right)^{h(x)}\left(f^{*}\right)^{h'(x)}$ we can easily calculate the higher order Bigeometric derivatives  (up to order three) as a product of Geometric multiplicative derivatives as 

\begin{eqnarray}
f^{*}(x)&=&\left(f^{\pi}(x)\right)^{\frac{1}{x}} \label{eq:fsfp}\\
f^{**}(x) &=&\left(\frac{f^{\pi\pi}(x)}{f^{\pi}(x)}\right)^{\frac{1}{x^{2}}}\label{eq:fsfp2} \\
f^{***}(x)&=&\left(\frac{f^{\pi\pi\pi}(x)\left(f^{\pi}(x)\right)^{2}}{\left(f^{\pi\pi}(x)\right)^{3}}\right)^{\frac{1}{x^3}}\label{eq:fsfp3}\\
\vdots &=& \vdots \nonumber
\end{eqnarray}

As indicated above we want to substitute the Geometric-multiplicative derivatives in the Multiplicative Taylor Theorem, given below, by its Bigeometric counterparts. The multiplicative Taylor theorem is given in \cite{BKO} as:

\begin{thm}[Multiplicative Taylor Theorem]
Let $A$ be an open interval
and let $f:A\rightarrow R$ be $n+1$ times {*} differentiable on
$A$. Then for any $x,x+h\in A$, there exists a number $\theta\in(0,1)$
such that 

\begin{equation}
f(x+h)=\prod_{m=0}^{n}\left(\left(f^{*(m)}(x)\right)^{\frac{h^{m}}{m!}}\right)\cdot\left(\left(f^{*(n+1)}(x+\theta h)\right)^{\frac{h^{n+1}}{(n+1)!}}\right)
\label{eq:mTaylor}
\end{equation}

\end{thm}

Equations  \eqref{eq:fsfp} -  \eqref{eq:fsfp3} suggest that the  $n$-th order Geometric-multiplicative derivative can be expressed in terms of the Bigeometric derivatives using the unsigned Stirling Numbers first kind $s(n,j)$  \cite{stirling} as stated in the following theorem.

\begin{thm}[Relation between Geometric Multiplicative and Bigeometric  derivative]
\label{thm:relgeobigeo}
The $n$-th Geometric multiplicative derivative can be expressed as a product of Bigeometric derivatives up to order $n$ as:
\begin{equation}
f^{*(n)}(x)=\left(\prod_{j=1}^{n}(f^{\pi(j)}(x))^{(-1)^{n-j}s(n,j)}\right)^{\frac{1}{x^{n}}}\label{eq:fsnfp}.
\end{equation}
\end{thm}

The proof of Theorem \ref{thm:relgeobigeo} can be simply carried out using mathematical induction as given in Appendix \ref{pgeobigeo}.

Finally, we substitute the Geometric multiplicative derivatives in the multiplicative Taylor theorem \eqref{eq:mTaylor} by the product of Bigeometric derivatives using \eqref{eq:fsnfp} we get:
\[
f(x+h)=\prod_{m=0}^{\infty}\left(\left(f^{*(m)}(x)\right)^{\frac{h^{m}}{m!}}\right)=\prod_{m=0}^\infty \left(\prod_{j=1}^m (f^{\pi (j)}(x))^{(-1)^{m-j}s(m,j)/x^m}\right)^{\frac{h^m}{m!}}
\]

Rearranging the factors in terms of the orders of the Bigeometric derivatives we get:
\begin{equation}
f(x+h)=\prod_{m=0}^\infty \left(\prod_{j=m}^\infty (f^{\pi (m)}(x))^{\frac{(-1)^{m-j}s(j,m)h^j}{x^j j!}}\right) = \prod_{m=0}^\infty \left((f^{\pi (m)}(x))^{\sum_{j=m}^\infty \frac{(-1)^{m-j}s(j,m)h^j}{x^j j!}}\right) 
. \label{eq:vTaylorx}
\end{equation}

With 
\begin{equation}
\label{eq:stirling2}
\sum_{j=m}^\infty (-1)^{j-m} s(j,m) \frac{x^j}{j!} = \frac{(\ln(1+x))^m}{m!},
\end{equation}
equation \eqref{eq:vTaylorx} simplifies to 
\begin{equation}
f(x+h)= \prod_{m=0}^\infty \left((f^{\pi (m)}(x))^{\frac{\left(\ln\left(1+\frac{h}{x}\right)\right)^m}{m!}}\right).
\end{equation}

Finally we can summarise the Bigeometric Taylor theorem as following. 

\begin{equation}
f(x+h)=\prod_{i=0}^{\infty}\left(f^{\pi(i)}(x)\right)^{\frac{(\ln(1+\frac{h}{x}))^{i}}{i!}}
\end{equation}

\cite{A} determined the Bigeometric Taylor theorem up to order 5 in $h/x$. Expansion of the logarithms up to order 5 in $h/x$ resembles the result of \cite{A}.

\begin{thm}[Bigeometric Taylor Theorem]
Let $A$ be an open interval
and let $f:A\rightarrow R$ be $n+1$ times ${\pi}$ differentiable on
$A$. Then for any $x,x+h\in A$, there exists a number $\theta\in(0,1)$
such that 

\begin{equation}
f(x+h)=\prod_{i=0}^{n}\left(f^{\pi(i)}(x)\right)^{\frac{(\ln(1+\frac{h}{x}))^{i}}{i!}}\cdot\left(\left(f^{\pi(n+1)}(x+\theta h)\right)^{\frac{(\ln(1+\frac{h}{x}))^{n+1}}{(n+1)!}}\right)
\label{eq:vTaylor}
\end{equation}

\end{thm}

While the Taylor expansion in \cite{A} of $f(x(1+\epsilon))$ is for small $\epsilon$, the Taylor expansion \eqref{eq:vTaylor} of $f(x+h)$ is for small $h$. Obviously, the Taylor expansion in \cite{A} can be written analogously as following:

\begin{equation}
f(x(1+\epsilon)) = \prod_{i=0}^{n}\left(f^{\pi(i)}(x)\right)^{\frac{(\ln(1+\epsilon)^{i}}{i!}}\cdot\left(\left(f^{\pi(n+1)}(x(1+\theta \epsilon))\right)^{\frac{(\ln(1+\epsilon)^{n+1}}{(n+1)!}}\right)
\end{equation}

\section{Bigeometric Runge-Kutta Method}
\label{sec3}

In this section we will present the Bigeometric Runge-Kutta method as an application of the Bigeometric Taylor theorem. The Newtonian Runge-Kutta Method is a widely used method for the numerical solution of initial value problems. In analogy to the Runge-Kutta Method in the framework of Newtonian Calculus and the Geometric multiplicative Runge-Kutta-Method \cite{RA13}, we will derive the Bigeometric Runge-Kutta Method in the framework of Bigeometric Calculus. Our starting point is here, as in the ordinary case, a Bigeometric initial value problem
\begin{equation}
y^\pi(x) = f(x,y) , \text{ with the initial value } y(x_0) = y_0.
\label{eq:vininitial}
\end{equation}

In order to explain the idea of the derivation of the Bigeometric Runge-Kutta method, the derivation will be carried out explicitly in the case of the  Bigeometric Euler Method (2nd order Bigeometric Runge-Kutta Method), and as the derivation is carried out analogously, only a summary of the derivation of the 4th order Bigeometric Runge-Kutta method will be given.

\subsection{2nd order Bigeometric Runge-Kutta Method or Bigeometric Euler Method }

In the ordinary case, the easiest approach to find an approximation to the solution of \eqref{eq:vininitial}, is called Euler Method. The Bigeometric Euler method will be derived in the following explicitly in analogy to the Euler method in ordinary calculus. For the Bigeometric Euler Method let us make the following ansatz:

\begin{equation}
y(x+h)=y(x)\cdot f_{0}^{a\ln(1+\frac{h}{x})}f_{1}^{b\ln(1+\frac{h}{x})}
\label{eq:eulera}
\end{equation}
with
\begin{eqnarray}
f_{0}&=&f(x,y), \\
f_{1}&=&f(x+ph,yf_{0}^{\frac{qh}{x}}).
\label{eq:f1}
\end{eqnarray}

The Bigeometric Taylor expansion \eqref{eq:vTaylor} for $y(x+h)$ up to order 2  is given as

\begin{equation}
y(x+h)=y(x)\left(y^{\pi}(x)\right)^{\ln(1+\frac{h}{x})}\left(y^{\pi\pi}(x)\right)^{\frac{1}{2!}\left[\ln(1+\frac{h}{x})\right]^{2}}  .
\label{eq:yvTaylor1}
\end{equation}

Substituting  \eqref{eq:vininitial} into \eqref{eq:yvTaylor1}, the Bigeometric Taylor expansion \eqref{eq:yvTaylor1} becomes:
\begin{equation}
y(x+h)=y(x)\left(f(x,y)\right)^{\ln(1+\frac{h}{x})}\left(f^{\pi}(x,y)(x)\right)^{\frac{1}{2!}\left[\ln(1+\frac{h}{x})\right]^{2}}  . 
\label{eq:yvTaylor2}
\end{equation}

In order to be able to compare  the equations \eqref{eq:yvTaylor2} and \eqref{eq:eulera} we have to expand $f_1$ from \eqref{eq:f1} also and substitute the result into \eqref{eq:yvTaylor2}. 

Using Bigeometric Taylor Theorem we obtain the expansion for $f_{1}$ 

\begin{equation}
f_{1}=f(x,y)\left[f_{x}^{\pi}(x,y)^{p}.f_{y}^{\pi}(x,y)^{y\frac{q}{x}\ln(f_{0})}\right]^{\ln(1+\frac{h}{x})}.
\label{eq:vtf1}
\end{equation}


Substituting \eqref{eq:vtf1} into \eqref{eq:eulera} the ansatz for the Bigeometric Euler Method becomes

\begin{equation}
y(x+h)=y(x)\cdot f(x,y)^{(a+b)\ln(1+\frac{h}{x})}\cdot f_{x}^{\pi}(x,y)^{bp\left(\ln(1+\frac{h}{x})\right)^{2}}\cdot f_{y}^{\pi}(x,y)^{bq\frac{y}{x}\ln(f(x,y)\left(\ln(1+\frac{h}{x})\right)^{2}}. 
\label{eq:eulerv2}
\end{equation}

Comparing the powers from the Bigeometric Taylor theorem \eqref{eq:yvTaylor2} with the powers from the Bigeometric Euler ansatz we can easily identify the following relations:

\begin{eqnarray}
a+b&=&1 \label{eq:veu1}\\
pb&=&\frac{1}{2}\label{eq:veu2}\\
qb&=&\frac{1}{2} \label{eq:veu3}
\end{eqnarray}

As the number of unknowns is greater than the number of equations, obviously we have infinitely many solutions of the equations \eqref{eq:veu1}-\eqref{eq:veu3}. One possible selection of the parameters $a,b,p,$ and $q$ can be  
\[
a=b= \frac{1}{2}, \quad \text{and } p=q=1.
\]

The parameters can be chosen differently according to the nature of the problem to be solved.

\subsection{4th order Bigeometric Runge-Kutta Method}

In science and engineering, generally the 4th order Runge-Kutta method is preferred, because it gives the most accurate approximation to initial value problems with reasonable computational effort. Also in the case of the Bigeometric Runge Kutta method, the 4th order method turns out to show the best balance between accuracy and computational effort. Consequently,  in analogy to the 2nd order Bigeometric Runge-Kutta method the starting point is again the Bigeometric Taylor expansion of $y(x+h)$, in this case, up to order four in $\ln\left(1+\frac{h}{x}\right)$

\begin{equation}
y(x+h)=y(x)\cdot\left(f^{\pi}(x)\right)^{\ln(1+\frac{h}{x})}\cdot\left(f^{\pi\pi}(x)\right)^{\frac{1}{2!}\left[\ln(1+\frac{h}{x})\right]^{2}}\cdot\left(f^{\pi(3)}(x)\right)^{\frac{1}{3!}\left[\ln(1+\frac{h}{x})\right]^{3}}\cdot\left(f^{\pi(4)}(x)\right)^{\frac{1}{4!}\left[\ln(1+\frac{h}{x})\right]^{4}} \label{eq:voltexp4}
\end{equation}

In analogy to the 2nd order case the ansatz for the 4th order Bigeometric Runge-Kutta method is

\begin{equation}
y(x+h)=y(x)\cdot f_{0}^{a\ln(1+\frac{h}{x})}\cdot f_{1}^{b\ln(1+\frac{h}{x})}\cdot f_{2}^{c\ln(1+\frac{h}{x})}\cdot f_{3}^{d\ln(1+\frac{h}{x})}
\label{eq:vrk4ans}
\end{equation}

with $f_0, f_1,f_2,$ and $f_3$ defined as following:

\begin{eqnarray}
f_{0}&=&f(x,y) \\
f_{1}&=&f\left(x+ph,yf_{0}^{\frac{qh}{x}}\right) \\
f_{2}&=&f\left(x+p_{1}h,yf_{0}^{\frac{q_{1}h}{x}}f_{1}^{\frac{q_{2}h}{x}}\right) \\
f_{3}&=&f\left(x+p_{2}h,yf_{0}^{\frac{q_{3}h}{x}} f_{1}^{\frac{q_{4}h}{x}} f_{2}^{\frac{q_{5}h}{x}}\right) .
\end{eqnarray}

After expanding  $f_1$, $f_2$, and $f_3$ using the Bigeometric Taylor theorem and substituting this expansions into \eqref{eq:vrk4ans} we can compare the powers of the Bigeometric derivatives with the ones in \eqref{eq:voltexp4}, and get the following relationships for the parameters: 
\begin{eqnarray}
p&=&q\\
p_1&=& q_1+q_2\\
p_2&=& q_3+q_4+q_5
\end{eqnarray}
and \begin{eqnarray}
a+b+c+d&=&1 \label{m4e13}\\
bp+cp_1+dp_2&=&\frac{1}{2}  \label{m4e23}\\
bp^2+cp_1^2+dp_2^2&=&\frac{1}{3}  \label{m4e33}
\end{eqnarray} 

The solution of \eqref{m4e13}-\eqref{m4e33} for $b$, $c$, and $d$ as functions of  $a$, $p_1$, and $p_2$ can be easily represented using   the Bigeometric Butcher Tableau \cite{butcher75}
\[
\begin{array}{c|cccc}
0& & & &\\
p& q& & &\\
p_1& q_1& q_2& &\\
p_2& q_3& q_4& q_5&\\
\hline
& a & b& c & d 
\end{array}.
\]

Also in this case the number of equations is more than the number of unknowns, therefore we have infinitely many solutions to the equations above. A suitable choice of the parameters actually depend also on the nature of the problem. In the case of the ordinary Runge-Kutta method, the following set of parameters is widely used, i. e.    
\begin{eqnarray} 
a&=&d=\frac{1}{6},\label{eq:ad}\\
b&=&c=\frac{1}{3},\\
p&=&p_{1}=q=q_{2}=\frac{1}{2},\\
p_{2}&=&q_{5}=1, \text{ and }\\
q_{1}&=&q_{3}=q_{4}=0.\label{eq:q1q3q4}
\end{eqnarray}

 The function is evaluated at four positions, i.e. at $x$, $x+ph$, $x+p_1 h$, and $x+p_2 h$. Reasonably $p_2=1$ so that we evaluate the function at the beginning and the end of the interval $[x, x+h]$. We select $p=p_1=\frac{1}{2}$ to calculate the function also in the middle of the interval. The weights of the contributions of $f_0$, $f_1$, $f_2$, and $f_3$ are $a$, $b$, $c$, and $d$ respectively. As $a+b+c+d=1$, we give equal weights for the end points of the interval, and put the emphasis on midpoint of the interval and get therefore $a=d=\frac{1}{6}$ and $b=c=\frac{1}{3}$. Nevertheless, the parameters can be selected in the framework of the  Butcher tableau for any problem independently to find the optimal solution.  

In the following we will use the parameters stated in \eqref{eq:ad}- \eqref{eq:q1q3q4} for the calculation of the  examples in the following section. 

 %

\section{Applications of the Bigeometric Runge-Kutta method}
\label{sec4}

In this section we want to apply the Bigeometric Runge-Kutta method to three examples. The first two examples are academic examples, where the solutions to the initial value problems are available in closed form. In these examples we will compare the results of the Bigeometric Runge-Kutta method with the ordinary Runge-Kutta method and determine the errors explicitly. The third example is a real world  example where the ordinary Runge-Kutta method breaks down in certain situations, whereas the Bigeometric Runge-Kutta method still gives accurate results. 
 
\subsection{Basic example}

As a simple and straightforward application of the 4th order Bigeometric Runge-Kutta method,  also used by  \cite{A}, as a simple academic example, we will consider the solution for the ordinary initial value problem
\begin{equation}
y'(x) = 1-\frac{1}{x}, \quad y(1)=1.
\label{eq:ano}
\end{equation}
Obviously, the corresponding Bigeometric initial value problem is
\begin{equation}
y^{\pi}(x)=\exp\left(\frac{x-1}{y}\right), \quad y(1)=1.
\label{eq:anv}
\end{equation}

The exact solution of the initial value problem \eqref{eq:ano} and \eqref{eq:anv}is 
\[
y(x) = x-\ln x.
\]

We will check the difference between the ordinary and the Bigeometric Runge-Kutta method by comparing the results of the initial value problems \eqref{eq:ano} and \eqref{eq:anv} exemplarily for a step size of $h=0.5$ and $n=6$ points, using the same parameters \eqref{eq:ad} - \eqref{eq:q1q3q4}.  In table \ref{tvrk} we will present the numerical results and their relative errors compared to the exact  result for the Bigeometric Runge-Kutta method and the ordinary Runge-Kutta method respectively.

\begin{table}[H]
\begin{center}
\begin{tabular}{llllll}
$x$ & $y_{exact}$& $y_{BRK4}$ & $y_{RK4}$ & relative error & relative error \\
 & &  &  & BRK4 &  RK4\\
\hline
1 & 1& 1 &1& 0& 0\\
1.5 & 1.09453& 1.10029 & 1.2123 & 0.00525979& 0.107595\\
2 & 1.30685 & 1.31299& 1.48915& 0.00469842& 0.139496\\
2.5 & 1.58371 & 1.58865& 1.80683& 0.0031184& 0.140885\\
3 &  1.90139 & 1.90483& 2.15268 & 0.00181087& 0.132162\\
3.5 & 2.24724 & 2.24927& 2.51915& 0.000905399& 0.120997 \\
4 &  2.61371& 2.61451 & 2.90136& 0.000307448& 0.110058\\
\hline 
\end{tabular}
\caption{\label{tvrk}Comparison of the results of the 4th order Bigeometric Runge-Kutta Method with the results of the 4th order ordinary Runge- Kutta Method with the exact values and their relative errors.}
\end{center}
\end{table}

The comparison of the results presented  in the table \ref{tvrk}  shows that the Bigeometric Runge-Kutta method gives more accurate results compared to the ordinary Runge Kutta method for the same parameters. The reduction of the step size by  one order in magnitude, i.e. $h=0.05$ we get for $y_{BRK4}(4) = 2.613727$ with a relative error of  $8.190459\times 10^{-6}$; on the other hand we get in the case of the ordinary Runge-Kutta  $y_{RK4}(4)=2.6500733$ with a relative error of $0.013914204$. Obviously reducing the step size by one order in magnitude has a more significant impact on the error in the case of the Bigeometric Runge-Kutta method as the error reduces by two orders in magnitude, whereas the error for the ordinary Runge-Kutta method reduces by only one order in magnitude, which is in parallel to the reduction of the step size. 

\begin{figure}[ht]
\centerline{\includegraphics{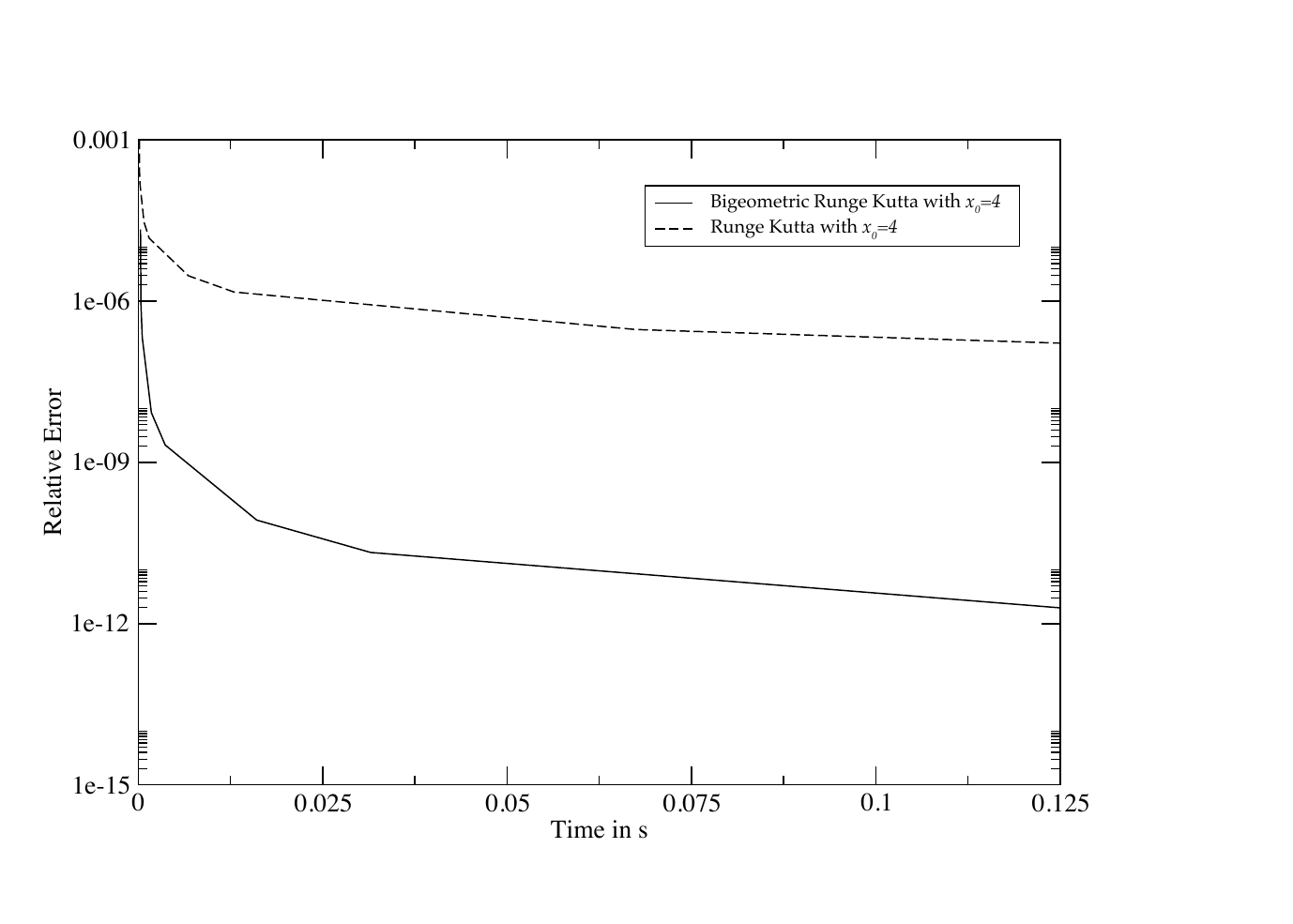}}
\caption{\label{aniszweska_perf_comp} Comparison of the  computation time vs. relative error for initial value $x_0=4$. of the Bigeometric Runge-Kutta method and the ordinary Runge-Kutta Method. (Intel Xeon 2.8GHz, g++ 4.81)}
\end{figure}

Obviously, in this example the  Bigeometric Runge-Kutta method achieves a higher accuracy compared to the ordinary Runge-Kutta method  by evaluating the function at fewer points. The basic operations of the Bigeometric Runge-Kutta method is multiplication, calculation of the exponential function, and calculation of the logarithm, compared to multiplication, addition, and subtraction in the ordinary case. Furthermore, we know that the multiplication of two $n$ bit numbers has the computational complexity  $O(n^2)$, addition and subtraction $O(n)$, and the calculation of the exponential function and logarithmic function of a $n$ bit input is $O(n^{5/2})$. Immediately  the question about the computation time arrises, or is the cost of using Bigeometric calculus more than reducing the value of $h$ accordingly. Therefore,  figure \ref{aniszweska_perf_comp} shows a comparison of the computation time vs. relative error for the initial value $x_0=4$. 

Evidently, one can observe that the relative error in the Bigeometric case is significantly smaller compared to the Newtonian case for the same computation time, which shows that at least in certain cases the Bigeometric Runge-Kutta method turns out to be a more efficient tool than the ordinary Runge-Kutta method.

\subsection{Non-Exponential Example}

Of course, one of the main concerns relating the previous example is that the solution involves a logarithm. So let us consider another example, where the solution will not be of exponential or logarithmic nature. 

Let us consider the following initial value problem

\begin{equation}
y'(x) =\frac{1}{2 y}, \quad y(4)=\sqrt{5}.
\label{eq:sqrt}
\end{equation}
and the corresponding Bigeometric initial value problem:
\begin{equation}
y^{\pi}(x)=\exp\left(\frac{x}{2 y^2}\right), \quad y(4)=\sqrt{5}.
\label{eq:sqrtpi}
\end{equation}

The solution of the initial value problems \eqref{eq:sqrt} and \eqref{eq:sqrtpi} can be easily determined as

\begin{equation}
y = \sqrt{1+x}.
\end{equation}

The numerical solutions of the Bigeometric initial value problem \eqref{eq:sqrtpi} and the ordinary initial value problem are summarised  in table \ref{sqrtvrk}.

\begin{table}[H]
\begin{center}
\begin{tabular}{llllll}
$x$ & $y_{exact}$ &$y_{VRK4}$ & $y_{RK4}$ & relative error  & relative error \\
 & & &  & BRK4 &  RK4\\
\hline
4 &2.23607& 2.23607 & 2.23607 & 0 & 0 \\
5 &  2.44949 & 2.44953 & 2.44037 & $1.75\times 10^{-5} $ & 0.00372\\
6 &  2.64575 & 2.64582  & 2.62947 & $2.58\times 10^{-5}$& 0.00615\\
7 & 2.82843 & 2.82851 & 2.80634 & $3.03\times 10^{-5} $ & 0.00781\\
8 &  3 & 3.0001 & 2.97307 & $3.29\times 10^{-5} $& 0.00898\\
9 &  3.16228 & 3.16239 & 3.13123 & $3.45 \times 10^{-5} $& 0.00982\\
10 &  3.31662 & 3.31674 & 3.28202 & $3.55 \times 10^{-5}$& 0.01043\\
\hline 
\end{tabular}
\caption{\label{sqrtvrk}Comparison of the results of the 4th order Bigeometric Runge-Kutta  and  Runge Kutta method with the exact values and its relative errors}
\end{center}
\end{table}

Obviously, we can see, that the relative error of the Bigeometric initial value problem is considerably smaller than in the ordinary case for the same step size $h$. As in the previous example, not the number of points or the value of $h$ alone determines the performance of the the method. Therefore, we will plot the computation time vs. relative error for both solutions for the initial value  $x_0=4$ and compare the numerical solutions, calculated with the corresponding Runge-Kutta method, and a fixed end point $x_{n+1}=5$ by varying the number of sub-intervals or $h$ and measure the computation time vs. the relative error. 

\begin{figure}[H]
\centerline{\includegraphics{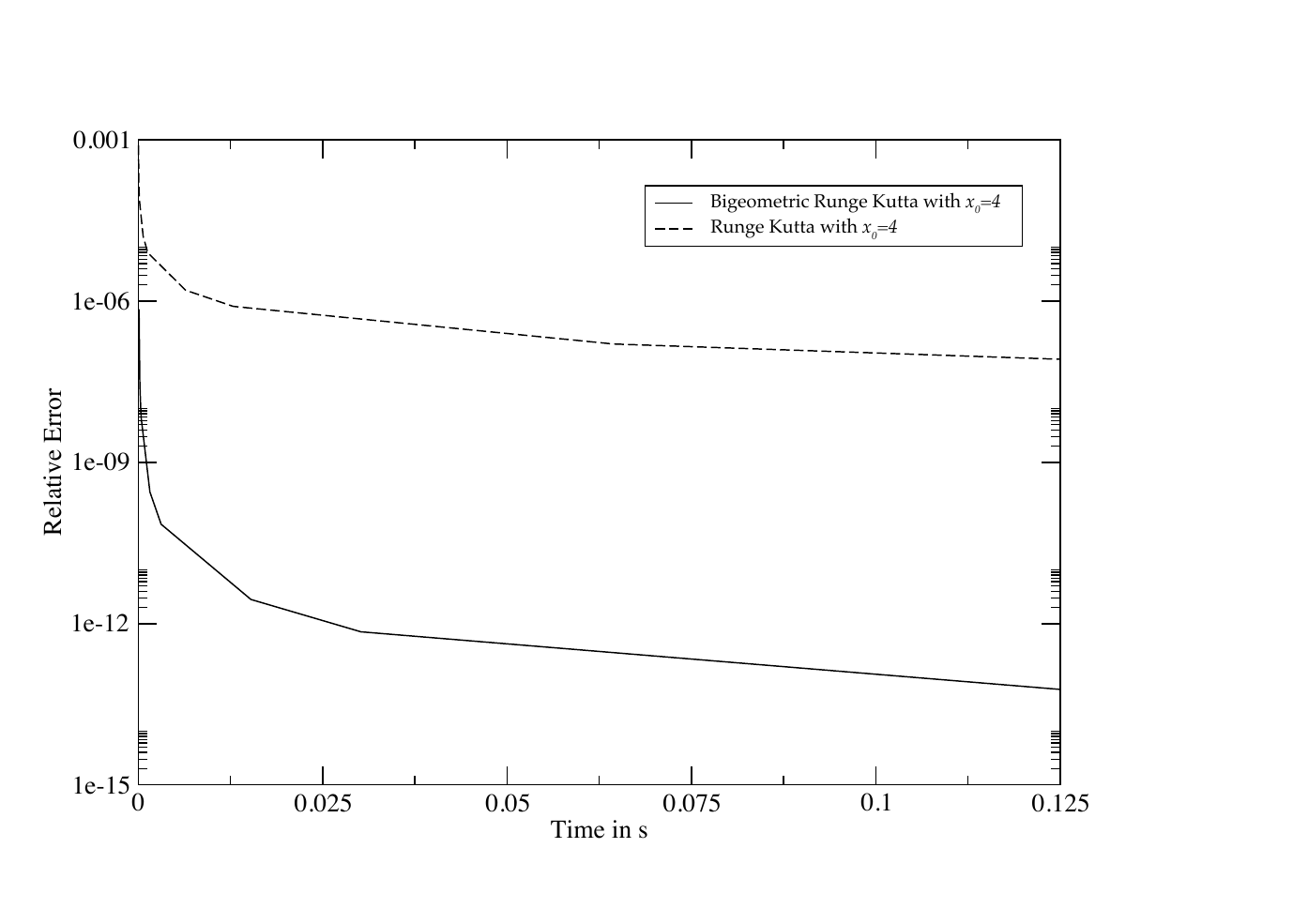}}
\caption{\label{sqrt}Comparison of the  computation time vs. relative error for initial value $x_0=4$. of the Bigeometric Runge-Kutta method and the ordinary Runge-Kutta Method. (Intel Xeon 2.8GHz, g++ 4.81)}
\end{figure}

Also in this case, the performance of the Bigeometric Runge-Kutta method seems to be superior to the one in the ordinary case. Of course there are also problems where the performance results are vice versa.

\subsection{Application to  Biological Modelling and its numerical results}

As an example for a coupled system of Bigeometric initial value problem, we will consider the mathematical model for tumor therapy with oncolytic viruses. 

\cite{agarwal11} developed a mathematical model of tumor therapy with oncolytic virus. The nonlinear model is based on a system of ordinary differential equations, modelling the size of the uninfected tumor cell population and the size of the infected tumor cell population. Agarwal \cite{agarwal11} carried out a stability analysis and checked their results using the  fourth order Runge-Kutta method. We will use the Bigeometric Runge-Kutta method and the ordinary Runge-Kutta method to calculate the size of the uninfected tumor cell population $x(t)$ and the size of the infected tumor cell population $y(t)$. Exemplarily we will carry out the comparison only for $y(t)$, as the results for  $x(t)$ are corresponding. The basic assumption in this model is that oncolytic viruses penetrate the tumor cells and replicate. Furthermore, infected tumor cells lead to infection of uninfected tumor cells with this oncolytic viruses. These oncolytic viruses preferably infect and lysis cancer cells and directly destruct them. In case of modification, anticancer proteins are produced. Based on these assumptions, the following set of ordinary differential equations are generate to model,  proposed by \cite{agarwal11}:

\begin{eqnarray}
\frac{dx}{dt}&=&r_{1}x\left(1-\frac{x+y}{K}\right)-\frac{bxy}{x+y+a}\label{eq:ago1}\\
\frac{dy}{dt}&=&r_{2}y\left(1-\frac{x+y}{K}\right)+\frac{bxy}{x+y+a}-\alpha y\label{eq:ago2}
\end{eqnarray}

with initial conditions: $x(0)=x_{0}>0$ and $y(0)=y_{0}>0$. First we have to clarify the parameters appearing in this nonlinear model.  $r_{1}$ and $r_{2}$ are maximum per capita growth rates of
uninfected and infected cells respectively. $K$ is the carrying capacity, $b$ is
the transmission rate, $a$ is the measure of the immune response
of the individual to the viruses which prevents it from destroying
the cancer and $\alpha$ is the rate of infected cell killing by the
viruses. All the parameters of the model are supposed to be nonnegative.

The corresponding system of Bigeometric differential equation system is given as

\begin{eqnarray}
x^{\pi}(t)&=&\exp\left[r_{1}t\left(1-\frac{x+y}{K}\right)-\frac{tby}{x+y+a}\right]\label{eq:agv1}\\
y^{\pi}(t)&=&\exp\left[r_{2}t\left(1-\frac{x+y}{K}\right)+\frac{tbx}{x+y+a}-\alpha t\right]\label{eq:agv2}
\end{eqnarray}

In order to find numerical approximations to  the functions $x(t)$ and $y(t)$ we used the ordinary and Bigeometric Runge Kutta methods and checked exemplarily for one set of parameters for what step size we get reasonable results. The time range is selected between $0$ and $1000$.  Because of the strongly nonlinear nature of the equations \eqref{eq:ago1}-\eqref{eq:ago2} and \eqref{eq:agv1}-\eqref{eq:agv2} we expect a small step size $h$ in both cases. Therefore, we carried out systematically gradual changes in the step size and the number of points to be calculated. 

\begin{figure}[H]
\centerline{\includegraphics[width=10cm]{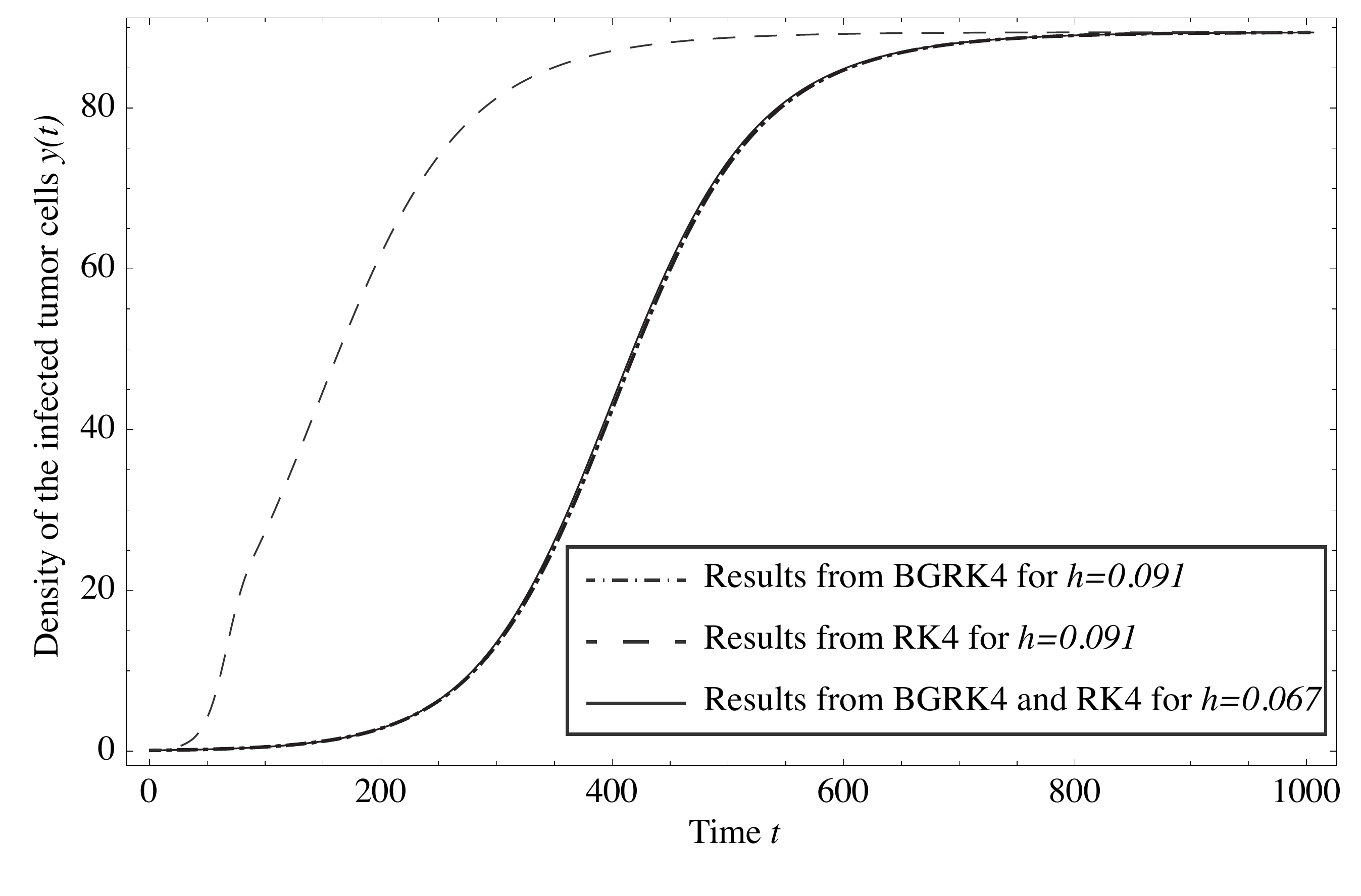}}
\caption{\label{infected} Density of infected tumor cells $y(t)$ as a function of time for the parameters $r_1=40$, $K=100$, $r_2=2$, $a=0.05$, $b=0.02$, and $\alpha=0.03$. The initial value $y(1)=0.1$. \label{figbigeo}}
\end{figure}

For a step-size $h=0.067$ the results from the Bigeometric Runge Kutta Method and the ordinary Runge-Kutta method coincide up to an absolute difference of $10^{-10}$. Therefore, as there is no closed form solution available for this problem, we accepted the solution for $h=0.067$ as the exact solution. With increasing step size up to $h=0.091$ the results of the Bigeometric Runge Kutta method differ not significantly from the ones for $h=0.067$, whereas the ordinary Runge Kutta method gives significantly different results as shown graphically in figure  \ref{figbigeo}. The absolute differences between the two computation experiments are summarised  in table \ref{tabbigeo}. Here, we can observe that the absolute difference for the Bigeometric Runge Kutta method is less than one, whereas the ordinary Runge Kutta method differs significantly more for the same step size $h=0.091$, i.e. up to 67.68.

\begin{table}[htdp]
\begin{center}
\begin{tabular}{llllll}
$t$ &  $y(t)$ & $y_{BG}(t)$ & $|y(t)-y_{BG}(t)|$& $y_{RK}(t)$ & $|y(t)-y_{RK}(t)|$ \\
\hline
1& 0.1& 0.1& 0& 0.1& 0\\   
100& 0.534879& 0.536288& 0.00140927& 27.2905& 26.7556\\   
200& 2.84556& 2.83341& 0.0121523& 61.8205& 58.9749\\   
300& 13.5502& 13.2125& 0.337775& 81.2303& 67.6801\\   
400& 43.4665& 42.5087& 0.957788& 87.0963& 43.6298\\   
500& 73.3533& 72.8296& 0.523669& 88.7399& 15.3866\\   
600& 84.8726& 84.7197& 0.152969& 89.2201& 4.34745\\   
700& 88.1091& 88.0662& 0.0428714& 89.3637& 1.25459
\end{tabular}
\end{center}
\caption{\label{infected} Density of infected tumor cells $y(t)$ as a function of time for the parameters $r_1=40$, $K=100$, $r_2=2$, $a=0.05$, $b=0.02$, and $\alpha=0.03$.The initial value $y(1)=0.1$.\label{tabbigeo}}
\label{tumor}
\end{table}

In order to compare the performance between the two methods, the absolute errors should be comparable. Therefore, we selected $h=0.0705$ for the ordinary Runge Kutta method and $h=0.091$ for the Bigeometric Runge Kutta method. As one can see from table \ref{bgrkrkperftab}, the maximum absolute difference is nearly twice of the absolute difference for the ordinary Runge Kutta method compared to the Bigeometric Runge Kutta method. On the other hand, the computation times are measured as $2.328$ seconds for the ordinary Runge Kutta method and $2.296$ seconds for the Bigeometric Runge Kutta method. So, even in this complicated mathematical model  the proposed method shows a higher performance at a higher accuracy.

\begin{table}[H]
\[
\begin{array}{llllll}
t & y(t)&y_{RK}(t) & |y(t)-y_{RK}(t)| & y_{BGRK}(t) &|y(t)-y_{BGRK}(t)|\\
\hline
1 & 0.1 &0.1 &0 & 0.1& 0\\
 100 & 0.53	4876 & 0.577645 & 0.0427696 & 0.536288 & 0.00141207 \\
 200 & 2.84555 & 3.09541 & 0.249859 & 2.83341 & 0.0121377 \\
 300 & 13.5502 & 14.5747 & 1.02456 & 13.2125 & 0.337715 \\
 400 & 43.4664 & 45.341 & 1.87463 & 42.5087 & 0.957674 \\
 500 & 73.3532 & 74.3369 & 0.983688 & 72.8296 & 0.523608 \\
 600 & 84.8726 & 85.1577 & 0.285053 & 84.7197 & 0.152951 \\
 700 & 88.1091 & 88.1891 & 0.0800294 & 88.0662 & 0.0428665 \\
 800 & 89.0339 & 89.0574 & 0.0234774 & 89.0213 & 0.0125593 \\
 900 & 89.3078 & 89.3148 & 0.00703348 & 89.304 & 0.00376082 \\
 1000 & 89.3901 & 89.3923 & 0.00212228 & 89.389 & 0.00113462 \\
\end{array}
\]
\caption{Comparison of the absolute errors of the results from the calculations of the ordinary Runge Kutta method with $h=0.0705$ ($n=14200$ points) and  the Bigeometric Runge Kutta Method for $h=0.091$ ($n=11000$ points). \label{bgrkrkperftab}}
\end{table}

Finally we can conclude that the Bigeometric Runge-Kutta method can be used to calculate approximate results for this model and that for a certain set of parameters the Bigeometric Runge-Kutta method gives better results than the ordinary Runge-Kutta method for larger step sizes. Of course, a more detailed analysis has to be carried out to show exactly for which class of initial value problems the Bigeometric Runge-Kutta method turns out to be more advantageous compared to the ordinary Runge-Kutta method. Nevertheless, we showed that there is a strong suspicion for certain problems the Bigeometric  Calculus can be a good base for the modelling and the numerical approximations of certain problems in science and engineering. 

The restriction to positive valued functions of real variable of the Bigeometric Calculus restricts the field of application drastically. Analogously to \cite{RA13} the theory of Bigeometric Calculus can be extended to complex valued functions of complex variable analogously to \cite{BR}. The problem that the Bigeometric derivative breaks down at the roots of the function, can be solved in analogy to \cite{RA13} by the application of the ordinary Runge-Kutta method on the interval $[\xi-h, \xi+h]$, if $\xi$ is a root of the function $f(x)$. With these two extensions the Bigeometric Runge-Kutta Method can be applied to any initial value problem.

\section{Conclusion}
In this paper we have stated and proven the differentiation rules for the Bigeometric derivative explicitly, and derived the Bigeometric Taylor theorem  on the basis of the Geometric multiplicative Taylor theorem exploiting the relation between the Geometric and Bigeometric multiplicative derivative. As an application of the Bigeometric Taylor expansion, we derived the Bigeometric Runge Kutta Method. The Bigeometric Runge-Kutta Method was tested on,  one hand for  basic examples and on the other, hand for the mathematical model of Agarwal for the Tumor Thereapy with Oncologic Virus \cite{agarwal11}. We chose two basic examples, one where the solution is of logarithmic nature and also used by  \cite{A}, second an example to show that the proposed method also works if the solution is not of logarithmic or exponential nature,  to check the performance of the proposed Bigeometric Runge-Kutta method. The comparison in the basic examples showed, that the Bigeometric Runge Kutta method gave better results for the same step size, and higher accuracy for a smaller computation time, which is an excellent indicator about the applicability of the proposed method.  In the case of the mathematical model of \cite{agarwal11} we could observe that the Bigeometric Runge-Kutta method gave better results for larger step sizes $h$ with a comparable computation time.  Finally we showed that the Bigeometric Runge-Kutta method is a serious tool for the solution of initial value problems, by also illuminating the computation time vs. relative error. 


\appendix
\section{Proofs of Bigeometric Differentiation Rules}
\label{bgdrproof}

\begin{enumerate}
\item Constant multiple Rule:
$$e^{x (\ln(cf(x)))^{'}}=e^{x\left(\frac{1}{cf(x)}\cdot(cf(x)')\right)}=e^{x \frac{f'(x)}{f(x)}}=f^{\pi}(x)$$

\item Product Rule: 
$$e^{x(\ln(f(x)g(x)))'}=e^{x \left(\frac{1}{f(x)g(x)} \cdot (f'(x)g(x)+f(x)g'(x))\right)}=e^{x\left(\frac{f'(x)g(x)}{f(x)g(x)}+\frac{f(x)g'(x)}{f(x)g(x)}\right)}=f^{\pi}(x)g^{\pi}(x)$$

\item  Quotient Rule:
$$e^{x(\ln(\frac{f(x)}{g(x)}))^{'}}=e^{x\left(\frac{g(x)}{f(x)}\left(\frac{f'(x)g(x)-g'(x)f(x)}{g(x)^{2}}\right)\right)}=e^{x\frac{f'(x)}{f(x)}}e^{-x\frac{g'(x)}{g(x)}}=\frac{f^{\pi}(x)}{g^{\pi}(x)}$$

\item Power Rule:

\begin{multline*}
e^{x \left(\ln(f(x))^{h(x)}\right)'}=e^{x\left(\frac{1}{f(x)^{h(x)}}\left[e^{h(x) \ln f(x)}\right]\right)}=e^{x\left(\frac{1}{f(x)^{h(x)}}f(x)^{h(x)}\left(h'(x)\ln f(x)+h(x)\frac{f'(x)}{f(x)}\right)\right)}=\\=e^{x \left(h'(x)\ln f(x)+ h(x)\frac{f'(x)}{f(x)}\right)}=\left(f(x)^{x h'(x)}\right)\left(f^{\pi}(x)\right)^{h(x)}
\end{multline*}

\item  Chain Rule: 
$$e^{x \ln(f\circ h)'(x)}=e^{x \left(\frac{1}{f(h(x))}\left(f'(h(x))h'(x)\right)\right)}=f^{\pi}(h(x))^{h'(x)}$$

\item Sum Rule: 
$$e^{x(\ln(f(x)+g(x)))'}=e^{x \left(\frac{1}{f(x)+g(x)} (f'(x)+g'(x))\right)}=\left(f^{\pi}(x)\right)^{\frac{f(x)}{f(x)+g(x)}} \left(g^{\pi}(x)\right)^{\frac{g(x)}{f(x)+g(x)}}$$
\end{enumerate}

\section{Proof of relationship between Geometric and Bigeometric derivative}
\label{pgeobigeo}

\begin{proof}
First let us check the formula for the first non-trivial case $n=2$. 
\begin{equation}
f^{**}(x) = \left(f^{\pi}(x)^{(-1)^{2-1} s(2,1)} \cdot f^{\pi\pi}(x)^{(-1)^{2-2} s(2,2)}\right)^{1/x^2} = \left(\frac{f^{\pi\pi}(x)}{f^{\pi}(x)}\right)^{\frac{1}{x^{2}}}\label{eq:fsfp2s}
\end{equation}
Equations \eqref{eq:fsfp2} and \eqref{eq:fsfp2s} are obviously identical.

Let \eqref{eq:fsnfp} be true for $n$ and check if it is true for $n+1$.

\begin{eqnarray*}
f^{*(n+1)} (x)&=& \frac{d^*}{dx^*} f^{*(n)} (x)\\
&=& \frac{d^*}{dx^*}\left(\prod_{j=1}^{n}(f^{\pi(j)}(x))^{\frac{(-1)^{n-j}s(n,j)}{x^n}}\right)
\end{eqnarray*}

With 
\[
\frac{d^* f(x)}{dx^*} = \left(\frac{d\pi f(x)}{dx^\pi}\right)^{1/x}
\]
we can calculate the $\pi $-derivative of the product.
\begin{eqnarray*}
f^{*(n+1)} (x)&=& \left[\frac{d^\pi}{dx^\pi}\left(\prod_{j=1}^{n}(f^{\pi(j)}(x))^{\frac{(-1)^{n-j}s(n,j)}{x^n}}\right)\right]^{1/x}\\
&=& \left[\prod_{j=1}^{n}(f^{\pi(j+1)}(x))^{\frac{(-1)^{n-j}s(n,j)}{x^n}} \cdot (f^{\pi(j)}(x))^{\frac{x (-n) (-1)^{n-j}s(n,j)}{x^{n+1}}}\right]^{1/x}\\
&=&\left[\prod_{j=2}^{n+1}(f^{\pi(j)}(x))^{(-1)^{n+1-j}s(n,j-1)} \cdot \prod_{j=1}^{n}(f^{\pi(j)}(x))^{(-1)^{n+1-j} n s(n,j)}\right]^{1/x^{n+1}}\\
&=&\left[\prod_{j=2}^{n}(f^{\pi(j)}(x))^{(-1)^{n+1-j}s(n,j-1)} \cdot(f^{\pi(n+1)}(x))^{(-1)^{n+1-(n+1}s(n,n)}  \cdot  \right. \\ 
&& \left. \cdot(f^{\pi(1)}(x))^{(-1)^{n} n s(n,1)} \cdot \prod_{j=2}^{n}(f^{\pi(j)}(x))^{(-1)^{n+1-j} n s(n,j)}\right]^{1/x^{n+1}} \\
&=&\left[\prod_{j=2}^{n}(f^{\pi(j)}(x))^{(-1)^{n+1-j}(s(n,j-1)+n s(n,j))}\right.  \cdot\\
&&\left. \cdot(f^{\pi(n+1)}(x))^{(-1)^{n+1-(n+1}s(n,n)}  \cdot(f^{\pi(1)}(x))^{(-1)^{n} n s(n,1)} \right]^{1/x^{n+1}} 
\end{eqnarray*}

Using the recurrence relation for the unsigned Stirling Number of first kind \cite{stirling}
\begin{equation}
\label{eq:stirlingid}
s(n+1,j)= n s(n,j) + s(n,j-1),
\end{equation}
we can simplify
\[
s(n,j-1)+n s(n,j)= s(n+1,j),
\]
and
\[
n s(n,1) = s(n+1,1) - s(n,0) = s(n+1,1),
\]
as $s(n+1,0)=0$. Finally with $s(n,n)= s(n+1,n+1)$ we can simplify 
\[
f^{*(n+1)}(x) =\left[\prod_{j=1}^{n+1}(f^{\pi(j)}(x))^{(-1)^{n+1-j}s(n+1,j)} \right]^{1/x^{n+1}},
\]
which completes the proof.
\end{proof}
\bibliographystyle{plainnat}

\end{document}